\documentclass[11pt]{article}
\usepackage{graphicx}
\usepackage{latexsym}
\def\bega{\begin{array}}
\def\enda{\end{array}}

\def\begi{\begin{itemize}}
\def\endi{\end{itemize}}
\def\bfq{{\bf q}}

\def\bff{{\bf f}}

\def\forall{\hbox{for all}~}
\def\L{{\bf L}}
\def\ds{\displaystyle}

\def\argmax{\hbox{arg}\!\max}

\def\ve{\varepsilon}

\def\I{{\cal I}}

\def\implies{\Longrightarrow}

\def\diams{\diamondsuit}

\def\vs{\vskip 2em}

\def\v{\vskip 1em}
\def\O{{\cal O}}
\def\C{{\cal C}}

\def\bel{\begin{equation}\label}
\def\eeq{\end{equation}}
\def\sqr#1#2{\vbox{\hrule height .#2pt
\hbox{\vrule width .#2pt height #1pt \kern #1pt
\vrule width .#2pt}\hrule height .#2pt }}
\def\square{\sqr74}
\def\endproof{\hphantom{MM}\hfill\llap{$\square$}\goodbreak}

\begin{document}
\title{\bf The Riemann Solver for Traffic Flow at an Intersection 
with Buffer of Vanishing Size}

\author{Alberto Bressan$^{(*)}$ and Anders Nordli$^{(**)}$\\    \\
(*) Department of Mathematics, Penn State University,\\
University Park, Pa.~16802, U.S.A.\\ \, \\
(**) Department of Mathematical Sciences,\\ Norwegian University of Science and Technology,\\ NO-7491 Trondheim, Norway \\
\, \\
e-mails:~ bressan@math.psu.edu,~anders.nordli@math.ntnu.no}

\maketitle
\begin{abstract} The paper examines the model of traffic flow at an intersection
introduced in \cite{BN1}, containing a buffer with limited size.
As the size of the buffer approach zero, it is proved that 
the solution of the Riemann problem with buffer 
converges to a self-similar solution described by a specific Limit Riemann Solver
(LRS).  Remarkably, this new Riemann Solver depends 
Lipschitz continuously on all parameters.
\end{abstract}

\section{Introduction}
\label{sec:0}
\setcounter{equation}{0}
Starting with the seminal papers by Lighthill, Witham, and Richards \cite{LW, R},
 traffic flow on a single road has been modeled in terms of a scalar conservation law:
\bel{LWR}
\rho_{t}+(v(\rho)\rho)_x~=~0\,.\eeq
Here $\rho$ is the density of cars, while $v(\rho)$ is their velocity, 
which we assume depends 
of the density alone.
To describe traffic flow on an entire network of roads, 
one needs to further introduce a set of
boundary conditions at road junctions \cite{GP}.  
These conditions should relate the traffic densities on 
incoming roads $i\in\I$ and outgoing roads $j\in \O$, depending on two main parameters:
\begi
\item[(i)] Driver's turning preferences.
For every $i,j$, one should specify 
the fraction $\theta_{ij}\in [0,1]$ of drivers arriving to  
from the $i$-th road, who wish to turn into 
the $j$-th road.   
\item[(ii)] Relative priorities assigned to different incoming roads.  
If the intersection is congested,  these describe  
the maximum influx of cars arriving from  each road $i\in \I$, 
allowed to cross the intersection.
\endi

Various junction models of have been proposed in the literature
\cite{CGP, G, GP, HR}.   See also \cite{BCGHP} for a survey.  
A convenient approach is to 
introduce a {\em Riemann Solver}, i.e.~a rule that specifies how to construct 
the solution in the special case where the initial data is constant 
on every incoming and each outgoing road.   As shown in \cite{CGP}, 
as soon as  a Riemann Solver is given,  the general  Cauchy problem for traffic 
flow near a junction can be uniquely solved (under suitable assumptions).

The recent counterexamples in \cite{BY} show that, 
on a network of roads,  
in general the Cauchy problem can be ill posed. Indeed, two distinct
solutions can be constructed for the same measurable initial data.   
On a network with several nodes,  nonuniqueness can occur even if 
the initial data have small total variation.
To readdress this situation, in  
\cite{BN1} an alternative intersection model was proposed, introducing a buffer
of limited capacity at each road junction.  For this new model, 
given any $\L^\infty$ initial data, the Cauchy problem has a unique solution,
which is robust w.r.t.~perturbations of the data. Indeed, one has continuous  
dependence even w.r.t.~the topology of weak convergence.

A natural question, addressed in the present paper, is what happens in the limit
as the size of the buffer approaches zero. 
For Riemann initial data, constant along each incoming and outgoing road, we
show that this limit is described by a  {\it Limit Riemann Solver} (LRS)
which can  be explicitly determined.     See (\ref{bbf})--(\ref{of}) in Section~2.

We recall that, in a model without buffer, the initial conditions consist of 
the constant densities $\rho^\diams_k$ on all incoming and outgoing roads
$k\in \I\cup\O$, together with the drivers' turning preferences 
$\theta_{ij}^\diams$.
On the other hand, in the model with  buffer,
these initial conditions
comprise also the 
length of the queues $q_j^\diams$, $j\in\O$, 
inside the buffer.  One can think of $q^\diams_j$ as 
the number of cars already inside the intersection (say, a traffic circle)
at time $t=0$,
waiting to access the outgoing road $j$.
  Our main results can be summarized as follows. 
\begi 
\item[(i)] For any given Riemann data $\rho^\diams_k,\,\theta_{ij}^\diams$,
one can choose initial queue sizes $q_j^\diams$ such that, for all 
$t>0$ the solution of the problem with buffer is exactly the same as the 
solution determined by the Riemann Solver (LRS).
\item[(ii)]  For any Riemann data $\rho^\diams_k,\,\theta_{ij}^\diams$,
and any  initial queue sizes $q_j^\diams$, as $t\to\infty$ the 
solution of the problem with buffer approaches asymptotically the solution
determined by the Riemann Solver (LRS).
\endi
Using the fact that the conservation laws (\ref{LWR}) are invariant under
space and time rescalings, from (ii) we obtain a convergence result 
as the size  of the buffer approaches zero. 

Our present results apply only to solutions of the Riemann problem, i.e.~with 
traffic density which is initially constant along each road.  
Indeed, for a general Cauchy problem
the counterexamples in \cite{BY} remain valid also for the Riemann Solver (LRS),
showing that the initial-value problem with measurable initial 
data can be ill posed.   Hence no convergence result can be expected.
This should not appear as a paradox: for every positive size of the buffer, the Cauchy problem has a unique solution, depending continuously on the initial data.
However, as the size of the buffer approaches zero, the solution can become more and more sensitive to small changes in the initial conditions. In the limit, uniqueness is lost.

An extension of our results may be possible
in the case of initial data with bounded variation, 
for a network containing one single  node.   In view
of the results in \cite{CGP, GP}, 
we conjecture that in this case the solution to the Cauchy problem with buffer
converges to the solution determined by the Riemann Solver (LRS).
\v
\section{Statement of the main results}
\label{sec:1}
\setcounter{equation}{0}

\begin{figure}[htbp]
\centering
  \includegraphics[scale=0.40]{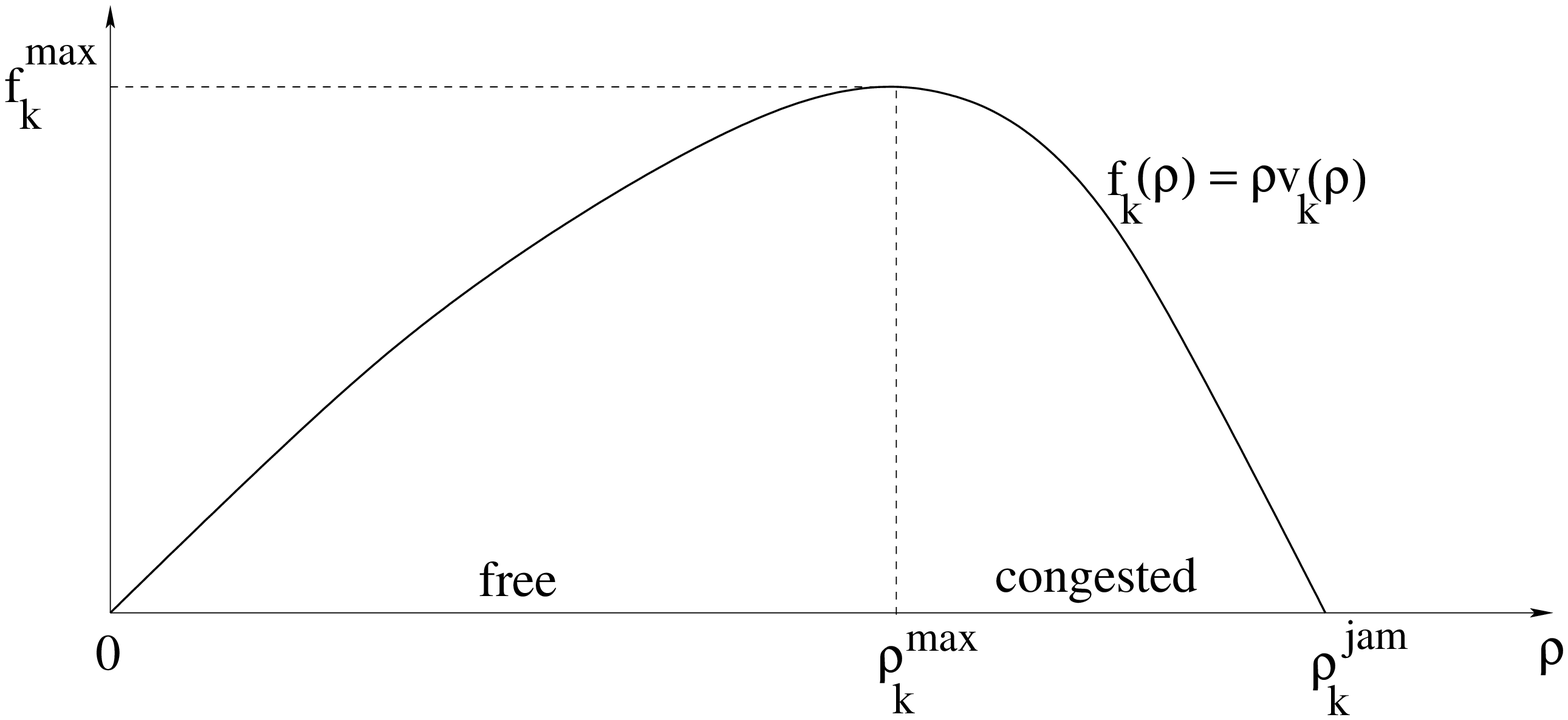}
    \caption{\small The flux $f_k$ as a function of the density $\rho$, along
    the $k$-th road.}
\label{f:tf34}
\end{figure}

Consider a family of $n+m$ roads, joining at a node.
Indices $i\in \{ 1,\ldots, m\}= \I$ denote {\it incoming roads}, 
while indices $j\in \{m+1,\ldots, m+n\}= \O$ denote {\it outgoing roads}.
On the $k$-th road, the density of cars $\rho_k(t,x)$ 
is governed by the scalar conservation law
\bel{clawi}\rho_t + f_k(\rho)_x~ =~ 0\,.\eeq
Here $t\geq 0$, while $x\in \,]-\infty,0]$ for incoming roads
and
$x\in [0,\,+\infty[\,$ for outgoing roads.
The flux function is $f_k(\rho) = \rho\,v_k(\rho)$, where $v_k(\rho)$ is the
speed of cars on the $k$-th road. 
We assume that 
each flux function $f_k$ satisfies
\bel{fi}
f_k\in \C^2, \qquad f_k(0)~= ~f_k(\rho_k^{jam}) ~=~0,\qquad f_k''(\rho)~<~0 \qquad
\hbox{for all}~~\rho\in [0, \rho_k^{jam}],\eeq
where $\rho_k^{jam}$ is the maximum possible density of cars on the $k$-th road.
Intuitively, this can be thought as a bumper-to-bumper packing, 
so that  the speed of cars is zero.
For a given road $k\in \{1,\ldots, m+n\}$, we denote by
$$f_k^{max}~\doteq~\max_s~ f_k(s)$$
the maximum flux and 
\bel{ri*}\rho_k^{max}~\doteq~\argmax_s ~f_k(s) \eeq
the traffic density corresponding to this maximum flux (see Fig.~\ref{f:tf34}).
 
Moreover, we say that 
$$\bega{rl}\rho~~\hbox{is a {\bf free} state if} &~~\rho\in \,[0, \, \rho_k^{max}]\,,\cr
 \rho~~\hbox{is a {\bf congested} state if} &~~\rho\in ~[\rho_k^{max}, \,\rho_k^{jam}]
\,.\enda$$



\v
Given initial data on each road
\bel{idi}
\rho_k(0,x) ~= ~\rho_k^\diams(x)\qquad\qquad k=1,\ldots,m+n,\eeq
in order to determine a unique solution to the Cauchy
problem we must supplement the conservation laws (\ref{clawi}) with
a suitable set of boundary conditions.  These provide additional constraints on 
the limiting values of the vehicle densities
\bel{brk}\bar \rho_k(t)~\doteq~\lim_{x\to 0} \rho_k(t,x)\qquad\qquad k=1,\ldots, m+n\eeq
near the intersection.
In a realistic model, these boundary conditions should depend on:
\begi

\item[{\bf (i)}] {\bf Relative priority given to incoming roads.}
For example, if the intersection is regulated by a crosslight, the flow  
will depend on the fraction $\eta_i\in \,]0,1[\,$ of time 
when cars arriving from the $i$-th road get a green light.

\item[{\bf (ii)}] {\bf Drivers' choices.} For every $i\in\I$, $j\in \O$, these are
modeled by assigning the fraction $\theta_{ij}\in [0,1]$ of drivers arriving 
from the $i$-th road who choose to
turn into the $j$-th road.  Obvious modeling considerations imply
\bel{Tij}
\theta_{ij}\in [0,1]\,,\qquad\qquad \sum_{j\in \O} \theta_{ij} ~=~1
\quad\hbox{for each}~i\in\I\,.\eeq
Since we are only interested in the Riemann problem, throughout the following
we shall assume that the $\theta_{ij}$ are given constants, satisfying (\ref{Tij}).
\endi



In \cite{BN1} a model of traffic flow at an intersection was introduced, 
including a buffer of limited capacity.   The incoming fluxes of cars 
toward the intersection are constrained by
the current degree of occupancy of the buffer.
More precisely, consider an intersection with $m$ incoming and $n$ outgoing roads.
The state of the buffer at the intersection is described by an $n$-vector
$$\bfq~=~(q_j)_{j\in\O}\,.$$
Here $q_j(t)$ is the number of cars at the intersection waiting to enter road $j\in\O$
(in other words, the length of the queue in front of road $j$).
Boundary values at the junction will be denoted by
\bel{traces}\left\{\bega{rl}\bar \theta_{ij}(t)&\doteq~\lim_{x\to 0-} \theta_{ij}(t,x),
\qquad i\in\I, ~j\in\O\,,\cr\cr
\bar \rho_i(t)&\doteq~\lim_{x\to 0-} \rho_i(t,x),\qquad i\in\I\,,\cr\cr
\bar \rho_j(t)&\doteq~\lim_{x\to 0+} \rho_j(t,x),\qquad j\in\O\,,\cr\cr
\bar f_i(t)&\doteq~f_i(\bar \rho_i(t))~=~\lim_{x\to 0-} f_i(\rho_i(t,x)),\qquad i\in\I\,,\cr\cr
\bar f_j(t)&\doteq~f_j(\bar \rho_j(t))~=~
\lim_{x\to 0+} f_j(\rho_j(t,x)),\qquad j\in\O\,.\enda\right.\eeq
Conservation of the total number of cars implies
\bel{BC}
\dot q_j(t)~=~\sum_{i\in\I} \bar f_i(t) \bar \theta_{ij} - 
\bar f_j(t)\qquad\qquad \forall j\in\O\,,\eeq
at a.e.~time $t\geq 0$.
Here and in the sequel, the upper dot denotes a derivative w.r.t.~time.
  Following \cite{GP}, we define the  maximum possible  flux
at the end of an incoming road as
\bel{omi}\omega_i~=~\omega_i(\bar \rho_i)~\doteq~\left\{ \bega{rl} f_i(\bar \rho_i) \quad&\hbox{if~$\bar \rho_{i}$
is a free state}, \cr
\cr
f_i^{max} \quad&\hbox{if~$\bar \rho_{i}$
is a congested state}, \enda\right.\qquad\quad i\in \I\,.\eeq  
  This is the largest flux $f_i(\rho)$ among 
all states $\rho$ that can be connected to $\bar \rho_i$ 
with a wave of negative speed.
Notice that the two right hand sides in (\ref{omi}) coincide if
$\bar \rho_i=\rho_i^{max}$.

Similarly, we define  the  maximum possible  flux
at the beginning of an outgoing road as 
\bel{omj}\omega_j~=~\omega_j(\bar \rho_j)~\doteq~
\left\{ \bega{rl} f_j(\bar\rho_{j}) \quad&\hbox{if~$\bar\rho_{j}$
is a congested state}, \cr
\cr
f_j^{max} \quad&\hbox{if~$\bar\rho_{j}$
is a free state}, \enda\right.\qquad\quad j\in\O\,.\eeq  


As in \cite{BN1}, we assume that the junction contains 
a buffer of size $M$.
 Incoming cars are admitted at a rate
depending of the amount of free space left in the buffer, regardless of their destination.
Once they are within the intersection, 
cars flow out at the maximum rate allowed by the outgoing road of their choice.
\v
{\bf Single Buffer Junction (SBJ).}  ~ {\it  
Consider a constant 
$M>0$, describing the maximum number of cars that can occupy the intersection
at any given time, and constants $c_i>0$, $i\in \I$,
accounting for priorities
given to different incoming roads.  

We then require that the  incoming fluxes $\bar f_i$ satisfy
\bel{bff}
\bar f_i~=~\min~\left\{ \omega_i\,,~~c_i\Big(M-\sum_{j\in\O} q_j\Big)\right\},
\qquad\qquad i\in\I\,.\eeq
In addition, the outgoing fluxes $\bar f_j$ should satisfy 
\bel{bqj}\left\{
\bega{l}\hbox{if $q_j>0$, then  $\bar f_j =\omega_j$,}\cr\cr
\hbox{if $q_j=0$, then $\bar f_j = 
\min\Big\{ \omega_j, ~\sum_{i\in \I} \bar f_i \bar\theta_{ij}\Big\}$,}
\enda\right.\qquad\qquad  j\in \O\,.\eeq
}
\v
Here $\omega_i=\omega_i(\bar\rho_i)$ and $\omega_j=\omega_j(\bar\rho_j)$
are the maximum fluxes defined at (\ref{omi})-(\ref{omj}).
Notice that {\bf (SBJ)} prescribes all the boundary fluxes $\bar f_k$, $k\in \I\cup\O$, depending on the boundary densities $\bar \rho_k$. 
It is natural to assume that the constants $c_i$ satisfy the  inequalities
\bel{Mi} c_i M~>~f_i^{max}\qquad\quad\forall i\in \I\,.\eeq
These conditions imply that, when the buffer is empty,  cars from
all incoming roads can access the intersection with the maximum possible flux
(\ref{omi}).  The analysis in \cite{BN1} shows that, with the above boundary
conditions, the Cauchy problem on a network of roads has a unique solution,
continuously depending on the initial data.

The main goal of this paper is to understand what happens when the 
size of the buffer approaches zero.   More precisely,  
assume that (\ref{bff})  is replaced by 
\bel{bfe}
\bar f_i~=~\min~\left\{ \omega_i\,,~~{c_i\over\ve}\Big(M\ve -\sum_{j\in\O} q_j\Big)\right\},
\qquad\qquad i\in\I\,.\eeq
Notice that (\ref{bfe}) models a buffer with size $M\ve$.   When $\sum_jq_j=M\ve$,
the buffer is full and no more cars are admitted to the intersection.

We will show that, as $\ve\to 0$, the limit of solutions to the Riemann problem with 
buffer of vanishing size
can be described by a specific Limit Riemann Solver.
\begi
\item[{\bf (LRS)}] {\it At time $t=0$, let the constant densities 
$\rho_i^\diams$, $\rho_j^\diams$ be given, together with 
drivers' preferences $\theta_{ij}$, ~ $i\in \I$, $j\in \O$.

Let $\omega_i^\diams=\omega_i(\rho_i^\diams)$ and  
$\omega_j^\diams=\omega_j(\rho_j^\diams)$  be the corresponding maximum
possible fluxes at the boundary of the incoming and outgoing roads,  as in
(\ref{omi})-(\ref{omj}).
Consider the one-parameter curve
$$s~\mapsto~\gamma(s) ~=~(\gamma_1(s),\ldots, \gamma_m(s)),$$
where 
$$\gamma_i(s)~\doteq~\min\{ c_i s\,,~ \omega_i^\diams\}.$$
Then for $t>0$ the Riemann problem is solved by the incoming fluxes
\bel{bbf} \bar f_i~=~\gamma_i(\bar s),\eeq
where
\bel{bars}\bar s~=~\max~\left\{ s\in [0,M]\,;~~\sum_{i\in \I} \gamma_i(s)\, \theta_{ij}~\leq~\omega_j^\diams
\quad\forall j\in\O\right\}.\eeq
In turn, by the conservation of the number of drivers, 
the outgoing fluxes are
\bel{of}
\bar  f_j~=~\sum_{i\in\I} \bar f_i\, \theta_{ij}\qquad\qquad j\in\O\,.\eeq
}
\endi
\v
By specifying all the incoming and outgoing fluxes $\bar f_i, \bar f_j$
at the intersection,
the entire solution of the Riemann problem is uniquely determined.
Indeed:

\begi
\item[(i)] For an incoming road $i\in \I$, there exists a unique boundary state 
$\rho_i^0=\rho_i(t, 0-)$ such that $f_i(\rho_i^0)= \bar f_i$ and moreover
\begi
\item[$\bullet$] If $\bar f_i=f_i(\rho_i^\diams)$, then $\rho_i^0=\rho_i^\diams$.
In this case the density of cars on the $i$-th road remains constant:
$\rho_i(t,x)\equiv \rho_i^\diams$. 

\item[$\bullet$] If $\bar f_i\not=f_i(\rho_i^\diams)$, then the solution to the
Riemann problem 
\bel{RPi}\rho_t+f_i(\rho)_x~=~0,\qquad\qquad \rho(0,x)~=~\left\{\bega{rl}
\rho_i^\diams&\qquad\hbox{if}~~x<0,\\[3mm]
\rho_i^0&\qquad\hbox{if}~~x>0,\enda\right.\eeq
contains only waves with negative speed.    In this case, the density
of cars on the $i$-th road coincides with the solution of (\ref{RPi}), 
for $x<0$.
\endi

\item[(ii)] For an outgoing road $j\in \O$, there exists a unique boundary state 
$\rho_j^0=\rho_j(t, 0-)$ such that $f_j(\rho_j^0)= \bar f_j$ and moreover
\begi
\item[$\bullet$] If $\bar f_j=f_j(\rho_j^\diams)$, then $\rho_j^0=\rho_j^\diams$.
In this case the density of cars on the $j$-th road remains constant:
$\rho_j(t,x)\equiv \rho_j^\diams$. 

\item[$\bullet$] If $\bar f_j\not=f_j(\rho_j^\diams)$, then the solution to the
Riemann problem 
\bel{RPj}\rho_t+f_j(\rho)_x~=~0,\qquad\qquad \rho(0,x)~=~\left\{\bega{rl}
\rho_j^0&\qquad\hbox{if}~~x<0,\\[3mm]
\rho_j^\diams&\qquad\hbox{if}~~x>0,\enda\right.\eeq
contains only waves with positive speed.    In this case, the density
of cars on the $i$-th road coincides with the solution of (\ref{RPj}),
for $x>0$.
\endi
\endi

\begin{figure}[htbp]
\centering
  \includegraphics[scale=0.5]{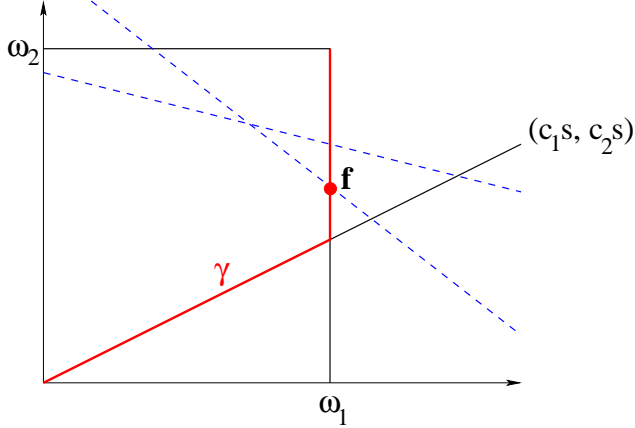}
    \caption{\small Constructing the solution of the the Riemann problem, according to the    limit Riemann solver (LRS),   
with two incoming and two outgoing roads.
    The vector $\bff=(\bar f_1,\bar f_2)$ of incoming fluxes
    is the largest point on the curve $\gamma$ that satisfies the two constraints
    $\sum_{i\in\I} \gamma_i(s) \theta_{ij}\leq \omega_j$, $j\in\O$. }
\label{f:tf145}
\end{figure}

{\bf Remark 1.} For the Riemann Solver constructed in \cite{BY}, the fluxes
$\bar f_k$ are {\em locally H\"older continuous} functions of 
the data $\rho_k^\diams, \theta_{ij}$,  on the domain where $\theta_{ij}>0$,
$\omega_j>0$ for all $j\in\O$.  \\ The Riemann Solver (LRS) has even better
regularity properties. Namely, the fluxes $\bar f_k$ defined at 
(\ref{bbf})--(\ref{of})
are {\em locally Lipschitz continuous}
functions of  $\rho_k^\diams, \theta_{ij}$.  Unfortunately, 
as remarked earlier, this additional regularity is still not sufficient
to guarantee the well-posedness of the Cauchy problem, for general
measurable initial data.
\v
Our first result refers to ``well prepared" initial data, where the initial lengths
of the queues are suitably chosen.
\v
{\bf Theorem 1.} {\it Let the assumptions  (\ref{fi}), (\ref{Mi}) hold.
Let Riemann data \bel{RD}\rho_k(0,x)= \rho_k^\diams\in [0, \rho_k^{jam}[\,,\qquad\qquad k\in \I\cup\O,\eeq
be assigned along each road, together
with drivers' turning preferences $\theta_{ij}$.

Then one can choose initial values $q_j^\diams$, $j\in\O$ for the queues 
inside the
buffer in such a way that the solution to the Riemann problem with buffer 
coincides with the self-similar solution determined 
by the Limit Riemann Solver (LRS).}
\v
Our second result covers the general case, where the initial sizes of the queues
are given arbitrarily, and the solution of the initial value problem 
with buffer is not self-similar.
\v
{\bf Theorem 2.}  {\it Let the assumptions  (\ref{fi}), (\ref{Mi}) hold.
Let Riemann data (\ref{RD})
be assigned along each road, together
with drivers' turning preferences $\theta_{ij}>0$ 
and initial queue sizes 
\bel{IQ}q_j(0)~=~q_j^\diams\,,\qquad\hbox{with} \quad\sum_{j\in \O}
q_j^\diams~<~M.\eeq
Then, as $t\to +\infty$, the solution $(\rho_k(t,x))_{k\in \I\cup\O}$ to the 
Riemann problem with buffer asymptotically converges to 
the self-similar solution $(\hat \rho_k(t,x))_{k\in \I\cup\O}$
determined by the Limit Riemann Solver (LRS).
More precisely:
\bel{lim}
\lim_{t\to +\infty}~{1\over t} \left(\sum_{i\in \I} \int_{-\infty}^0 
|\rho_i(t,x)-\hat \rho_i(t,x)|\, dx + \sum_{j\in \O} 
\int_0^{+\infty} |\rho_j(t,x)-\hat \rho_j(t,x)|\, dx\right)~=~0.\eeq
}
\v
A proof of the above theorems will be given in Sections 4 and 5, respectively.
By an asymptotic rescaling of time and space, using Theorem~2 we can describe
the behavior of the solution to a Riemann problem, 
as the size of the buffer approaches zero.
\v
{\bf Corollary (limit behavior for a buffer of vanishing size).} {\it Let $f_k, \theta_{ij}, c_i, M$  be as in Theorem~2. 
Let Riemann data (\ref{RD})
be assigned along each road, together
with drivers' turning preferences $\theta_{ij}>0$ 
and initial queue sizes as in (\ref{IQ}).

For $\ve>0$, let $(\rho^\ve_k(t,x))_{k\in \I\cup\O}$ be the solution
to the initial value problem 
with a buffer of  size $M\ve$,  obtained by replacing (\ref{bff}) with (\ref{bfe})
and choosing $q^\ve_j(0)=\ve q_j^\diams$ as initial sizes of the queues.

Calling  $\hat \rho_k$ the self-similar solution
determined by the Limit Riemann Solver (LRS) with the same initial data 
(\ref{RD}), for every $\tau> 0$ we have
\bel{lim2}
\lim_{\ve\to 0}~ \left(\sum_{i\in \I} \int_{-\infty}^0 
|\rho^\ve_i(\tau,x)-\hat \rho_i(\tau,x)|\, dx + \sum_{j\in \O} 
\int_0^{+\infty} |\rho^\ve_j(\tau,x)-\hat \rho_j(\tau,x)|\, dx\right)~=~0.\eeq
}

\v
{\bf Proof of the Corollary.}
Let $(\rho_k(t,x))_{k\in \I\cup\O}$ be the solution
constructed in Theorem~2.   Then, for every $i\in\I$, by a rescaling 
of coordinates we obtain 
$$\bega{l}\ds\lim_{\ve\to 0}~\int_{-\infty}^0 
|\rho^\ve_i(\tau,x)-\hat \rho_i(\tau,x)|\, dx ~=~\lim_{\ve\to 0}~
\int_{-\infty}^0 
\left|\rho_i\Big({\tau\over\ve},\,{x\over\ve}\Big)-\hat \rho_i\Big({\tau\over\ve},
\,{x\over\ve}\Big)
\right|\, dx\\[4mm]
\qquad\ds  =~\lim_{\ve\to 0}~\ve
\int_{-\infty}^0 
\left|\rho_i\Big({\tau\over\ve},\,x\Big)-\hat \rho_i\Big({\tau\over\ve},
\,x\Big)
\right|\, dx
=~\lim_{t\to\infty}~{\tau\over t}\,
\int_{-\infty}^0 
\bigl|\rho_i(t,\,x)-\hat \rho_i(t,
\,x)
\bigr|\, dx~=~0.
\enda$$
In the last step we used Theorem~2 in connection with
the variable change $t=\tau/\ve$.   For $j\in\O$,
the difference $|\rho^\ve_j-\hat \rho_j|$ is estimated in an entirely similar way.
\endproof

\v
\section{The Riemann problem with buffer}
\label{sec:2}
\setcounter{equation}{0}
We consider here an initial value problem with Riemann data,
so that the initial density is constant on every incoming and outgoing road.
\bel{RP1}\left\{\bega{rl}\rho_i(0,x)&=~\rho^\diams_i\qquad i\in\I\,,\cr\cr
\rho_j(0,x)&=~ \rho^\diams_j\,,\qquad j\in\O\,,\enda\right.
\qquad\qquad q_j(0) ~=~q^\diams_j\quad j\in\O\,.\eeq

\begin{figure}[htbp]
\centering
\includegraphics[scale=0.55]{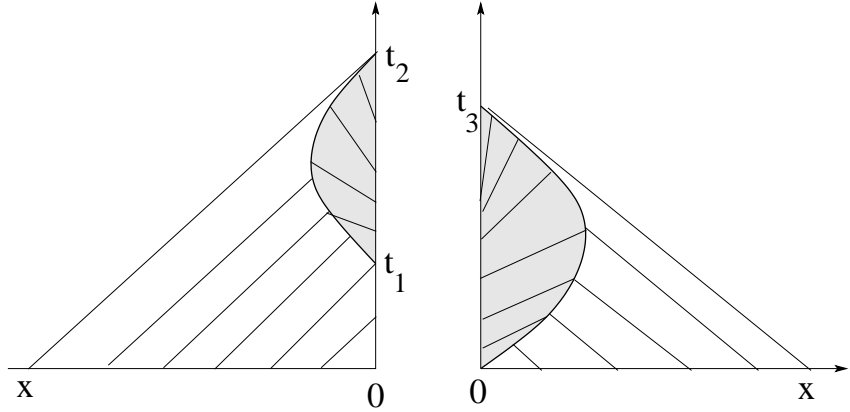}
    \caption{\small Left: an incoming road which is initially free.
    For $t_1<t<t_2$ 
    part of the road is congested (shaded area).
 Right: an outgoing road which is initially congested.
    For $0<t<t_3$ 
    part of the road is free (shaded area).  In both cases, a shock 
    marks the boundary
    between the free and the congested region.
}
\label{f:tf170}
\end{figure}

We decompose the sets of indices as
$$\I~=~\I^f\cup\I^c\,,\qquad\qquad \O~=~\O^f\cup\O^c,$$
depending on whether these roads are initially free or congested. 
More precisely:
\bel{FC}
\bega{rl}\I^f~\doteq~\{i\in \I\,;~~\rho_i^\diams < \rho_i^{max}\}\,,
\qquad\qquad &\O^f~\doteq~\{j\in \O\,;~~\rho_j^\diams \leq \rho_j^{max}\}\,,\cr\cr
\I^c~\doteq~\{i\in \I\,;~~\rho_i^\diams \geq \rho_i^{max}\}\,,\qquad
\qquad &
\O^c~\doteq~\{j\in \O\,;~~\rho_j^\diams > \rho_j^{max}\}\,.
\enda\eeq
Observe that
\begi
\item  If $i\in \I^c$, then the $i$-th incoming road will always be congested,
i.e.~$\rho_i(t,x)\geq \rho_i^{max}$ for all $t,x$.
\item  If $j\in \O^f$, then the $j$-th outgoing road will always be free,
i.e.~$\rho_j(t,x)\leq \rho_j^{max}$ for all $t,x$.
\item If $i\in \I^f$, then part of the $i$-th  road can become congested
 (Fig.~\ref{f:tf170}, left).
\item If $j\in \O^c$, then part of the $j$-th  road can become free
 (Fig.~\ref{f:tf170}, right).
\endi
 

\v
The next lemma plays a key role in the proof of Theorem 2.  It shows that, for any $t>0$,
 the maximum possible flux at  the boundary of any incoming 
or outgoing road is greater or equal to the maximum flux computed at $t=0$. 
\v
{\bf Lemma 1.} {\it Let $\rho_k=\rho_k(t,x)$, $k\in \I\cup\O$ be the solution 
of the Riemann problem with initial data (\ref{RP1}).   As in (\ref{omi})-(\ref{omj})
call $\omega_k^\diams= \omega_k(\rho_k^\diams)$ the maximum possible fluxes.  
 Similarly, for $t>0$ call $\omega_k(t)=\omega_k(\bar \rho_k(t))$ the corresponding
maximum fluxes. Then 
\bel{mf}
\omega_k(t)~\in~\bigl\{ \omega_k^\diams, \,f_k^{max}\bigr\}\qquad\qquad\forall k\in \I\cup\O,\quad t\geq 0.\eeq}
\v
{\bf Proof.} 
{\bf 1.} We first consider an incoming road $i\in\I$.

CASE 1: The road is initially congested, namely $\rho_i^\diams\geq \rho_i^{max}$.  
In this case the $i$-th 
road always remains congested and we have 
$\omega_i(t)= \omega_i^\diams = f_i^{max}$, for every $t\geq 0$.

CASE 2: The road is initially free, namely $\rho_i^\diams< \rho_i^{max}$.  
For a given $t>0$, two subcases may occur.   
\begi
\item[(i)] There exists a characteristic with positive speed, reaching the point $(t,0)$.
Since this characteristic must start at a point $x_0<0$, we conclude that 
$\rho_i(t, 0-) = \rho_i(0, x_0) = \rho_i^\diams$.
Hence $\omega_i(t)= \omega_i^\diams$.
\item[(ii)] There exists a  neighborhood of $(t,0)$ covered with characteristics having negative speed.   In this case $\rho_i(t, 0-) \geq \rho_i^{max}$, hence
$\omega_i(t) = f_i^{max}$.
\endi
\v
{\bf 2.} For an outgoing road $j\in \O$, the analysis is similar.

CASE 1: The road is initially free, namely $\rho_j^\diams\leq \rho_j^{max}$.  
In this case the $j$-th 
road always remains free and we have 
$\omega_j(t)= \omega_j^\diams = f_j^{max}$, for every $t\geq 0$.

CASE 2: The road is initially congested, namely $\rho_j^\diams>\rho_j^{max}$.  
For a given $t>0$, two subcases may occur.   
\begi
\item[(i)] There exists a characteristic with negative speed, reaching the point $(t,0)$.
Since this characteristic must start at a point $x_0>0$, we conclude that 
$\rho_j(t, 0+) = \rho_j(0, x_0) = \rho_j^\diams$.
Hence $\omega_j(t)= \omega_j^\diams$.
\item[(ii)] There exists a  neighborhood of $(t,0)$ covered with characteristics having
positive speed.   In this case $\rho_j(t, 0+) \geq \rho_j^{max}$, hence
$\omega_j(t) = f_j^{max}$.
\endi
\endproof

\section{Proof of Theorem 1}
\label{sec:4}
\begin{figure}[htbp]
\centering
\includegraphics[scale=0.5]{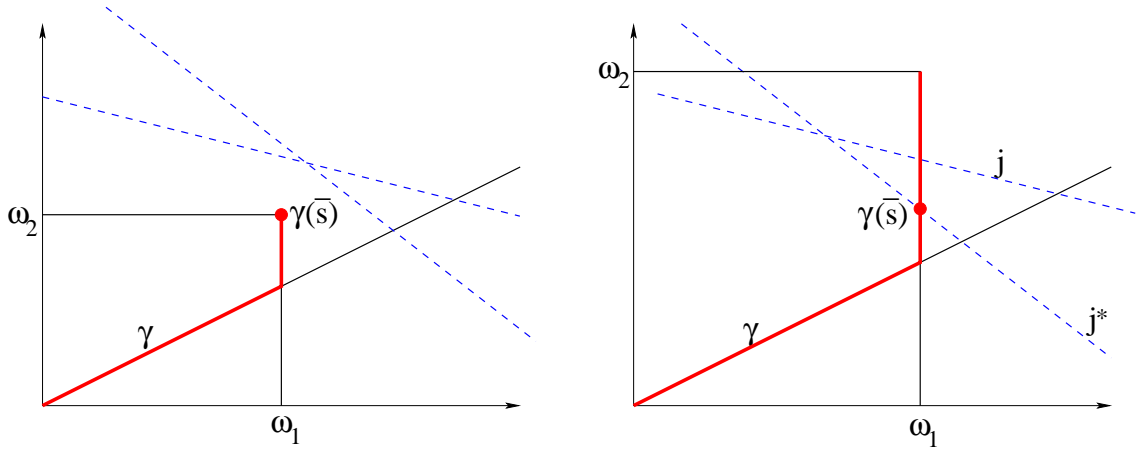}
   \caption{\small The two cases in the proof of Theorem~1.
   Left: none of the outgoing roads provides a restriction on the
   fluxes of the incoming roads. The queues are zero.   Right: one of the outgoing roads is congested and restricts the maximum flux through the node.}
\label{f:tf205}
\end{figure}

\setcounter{equation}{0}
Let  $\rho_k^\diams$, $k\in \I\cup\O$ be the initial densities of cars 
on the incoming and outgoing roads, and let $\theta_{ij}$ be the 
drivers' turning preferences, as in (\ref{Tij}). Call $\omega_i^\diams, 
\omega_j^\diams$
the maximum possible boundary fluxes on the incoming and outgoing roads,
and define $\bar s$ as in (\ref{bars}). Two cases will be considered, shown in Fig.~\ref{f:tf205}.
\v
CASE 1: $\bar s=M$, so that $\gamma(\bar s) ~=~
(\omega_1^\diams,\omega_2^\diams,\ldots, \omega_m^\diams)$.
This is the case where none of the incoming roads remains congested, 
and all the 
drivers arriving at the intersection can immediately proceed 
to the outgoing road of their choice.

In this case we choose the initial queues  
$$q^\diams_j~=~0\qquad\qquad \forall ~~ j\in \O\,.$$
With these choices,
the  solution of the Cauchy problem with buffer coincides with the self-similar 
solution determined by the Limit Riemann Solver (LRS).
The buffer remains always empty:
$q_j(t)=0$ for all $t\geq 0$ and $j\in\O$.
\v
CASE 2: $\bar s<M$.  In this case there exists an index $j^*\in\O$ such that 
\bel{gj*}\sum_{i\in \I} 
\gamma_i(\bar s)\theta_{ij^*} ~=~ \omega_{j^*}^\diams\,.\eeq
When this happens, the entire flow through the intersection is
restricted by the number of cars that can exit toward the single 
congested road $j^*$.  
We then define
\bel{q*}q^*~\doteq~M-\bar s\,,\eeq
and choose 
the initial queues to be
\bel{qu2}q_j^\diams ~=~ \left\{ \bega{cl} q^*\qquad &\hbox{if}~~~j=j^*\cr
0\qquad &\hbox{if}~~~j\not=j^*.\enda \right.\eeq
Then the corresponding solution coincides with the self-similar 
solution determined by the Limit Riemann Solver (LRS).
Indeed, by the definition of $\gamma(\bar s)$, for every $j\in \O$ we have
\bel{gi}\sum_i \min\Big\{ c_i(M-q^*), ~\omega_i\Big\}\cdot \theta_{ij}
~=~\sum_i\gamma_i(\bar s)\, \theta_{ij}~\leq~\omega_j\,,\eeq
with equality holding when $j=j^*$.
By (\ref{gi}),
all queues remain constant in time, namely
$q_{j^*}(t)= q^*$ and $q_j(t)=0$ for  $j\not= j^*$.

\endproof

\v
{\bf Remark 2.}  In the proof of Theorem~1,  the queue sizes $q_j^\diams$ 
may not be uniquely 
determined.  Indeed, in Case~2 
there may exist two distinct indices $j^*_1, j^*_2\in\O$ such that
$$\sum_{i\in \I} 
\gamma_i(\bar s)\theta_{ij_1^*} ~=~ \omega_{j_1^*}\,,\qquad\qquad 
\sum_{i\in \I} 
\gamma_i(\bar s)\theta_{ij_2^*} ~=~ \omega_{j_2^*}\,.$$ 
When this happens, we can choose the queue sizes to be
\bel{qu3}q_j^\diams ~=~ \left\{ \bega{cl} \alpha q^*\qquad &\hbox{if}~~~j=j^*_1\,,
\\[3mm]
(1-\alpha) q^*\qquad &\hbox{if}~~~j=j^*_2\,,\\[3mm]
0\qquad &\hbox{if}~~~j\notin\{j^*_1,j^*_2\},\enda \right.\eeq for any 
choice of $\alpha\in [0,1]$.
\v
\section{Proof of Theorem 2}
\label{sec:5}
\setcounter{equation}{0}

In this section we prove that, for any initial data,
as $t\to +\infty$  the solution to the Riemann
problem with buffer converges as to the self-similar
function determined by the Limit Riemann Solver (LRS).
The main argument can be divided in three main steps. (i) Establish 
an upper bound on the size $q=\sum_j q_j$ of the queue insider the buffer,
showing that $\limsup_{t\to\infty} q(t)\leq M-\bar s$.
(ii) Establish the lower bound $\liminf_{t\to\infty} q(t)\geq M-\bar s$.
(iii) Using the previous steps, show that as $t\to \infty$ all boundary
fluxes in the solution with buffer 
converge to the corresponding fluxes determined 
by (LRS).   From this fact, the limit (\ref{lim}) follows easily.
\v
Given the densities $\rho_i^\diams$ on the incoming roads   $i\in\I$,
call $\omega^\diams_i$ the 
corresponding maximal flows, as in (\ref{omi}).  
Call $\hat q_i$ the value of the queue such that 
$$c_i(M-\hat q_i)~ =~ \omega_i^\diams\,.$$
Without loss of generality, we can assume
\bel{hqi}0~\leq~\hat q_m~\leq~\cdots~\leq~\hat q_2~\leq~\hat q_1\,.\eeq

At an intuitive level, we have
\begi
\item  If the queue inside the buffer is small, i.e.~$q<\hat q_i$, then  all drivers
arriving from the $i$-th road can access the intersection, and the $i$-th road 
will become  free.
\item If the queue inside the buffer is large, i.e.~$q>\hat q_i$, 
then not all drivers coming from the $i$-th road can immediately access the intersection, 
and the $i$-th road will become congested.
\endi
This can be formulated in a more precise way as follows.
By the definition (\ref{bars}), if  $q> M-\bar s $
one has
\bel{in6}
\sum_{i\in \I} \min\{ c_i(M- q), ~\omega_i^\diams\}\cdot
\theta_{ij}~<~\omega_j^\diams\qquad \quad
\hbox{for every}~~j\in \O \,.\eeq
On the other hand, if $q<M-\bar s$, let $j^*\in\O$ be an index such that
(\ref{gj*}) holds.   We then have
\bel{in7}
\sum_{i\in \I} \min\{ c_i(M- q), ~\omega_i^\diams\}\cdot
\theta_{ij^*}~>~\omega_{j^*}^\diams\,.\eeq\v

\begin{figure}[htbp]
\centering
  \includegraphics[scale=0.5]{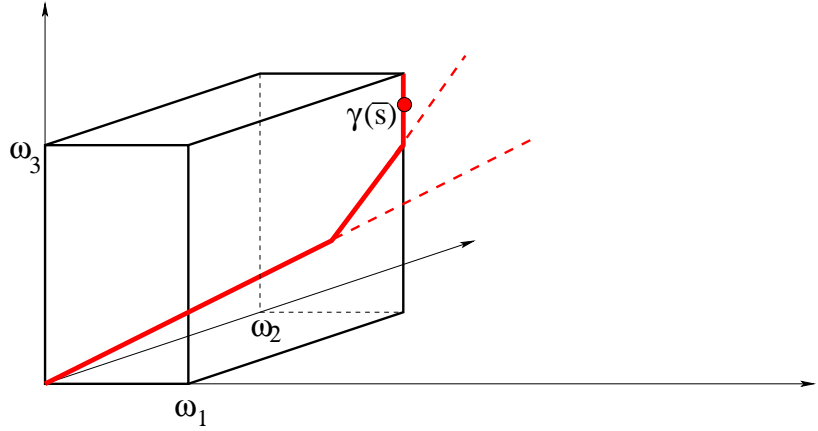}
    \caption{\small  A case with three incoming roads. For large times, the first two roads become free, while the third road remains congested. }
\label{f:tf207}
\end{figure}

The proof is achieved in several steps.
\v
{\bf 1.} We first study the case where, in the solution determined by the Limit Riemann Solver,  at least one of the 
outgoing roads is congested
(Fig.~\ref{f:tf205}, right), so that (\ref{gj*}) holds.  
Let $\bar s$ be as in (\ref{bars}). As in (\ref{q*}), 
define the 
asymptotic size of the queue to be $q^*= M-\bar s >0$.
To fix the ideas, assume 
\bel{hqel}0~\leq~\hat q_m~\leq~\cdots~\leq~
\hat q_{\nu+1}~\leq~ q^*~<~\hat q_\nu~\leq ~\cdots~\leq~\hat q_2~\leq~\hat q_1\,.\eeq
In this setting, we will show that for $t$ large the incoming roads
$i=1,\ldots, \nu$ will be free, while the 
incoming roads with $\hat q_i<q^*$ will be congested.
More precisely, we shall prove the following

{\bf Claim.} {\it 
There exist times
\bel{ttt}
0~= ~t_0~=~\tau_0~<~t_1~<~\tau_2~<~t_2~<~\cdots~<\tau_{\nu}~<~t_\nu\eeq
and constants $\delta_\ell,\ve_\ell>0$, $\ell=1,\ldots, \nu$, with the following properties.

\begi
\item[(i)] If $t\geq t_{\ell-1}$, then we have the implication
\bel{impk}
q(t)~\geq~ \hat q_\ell-\delta_\ell\qquad\implies\qquad
\dot q(t)~\leq ~-\ve_\ell~<~0\,.\eeq   

\item[(ii)]  If $t\geq \tau_\ell$, then $q(t)\leq \hat q_\ell-\delta_\ell$

\item[(iii)]   For all times $t\geq t_\ell$ the incoming road $\ell$ is free.  Hence its flux 
near the intersection satisfies
\bel{lfre}
f_\ell(t)~=~\omega_\ell^\diams\qquad\qquad\forall t\geq t_\ell\,.\eeq
\endi
}

The above claim is proved by induction on $\ell=1,\ldots,\nu$.

We begin with $\ell=1$.    For any $t\geq  0$, if $q(t)\geq \hat q_1$
then by (\ref{hqel}) we have $q(t)>q^*$.
Using Lemma~1, we thus obtain
$$\dot q_j(t) ~\leq~\sum_i c_i(M-q(t))\theta_{ij}  - \omega_j^\diams\qquad\qquad \hbox{if}~~q_j(t)>0.$$
Therefore, if $q_j(t)>0$, then 
$$\dot q_j(t) ~\leq~-2\ve_{1j}~<~0$$
for some $\ve_{1j}>0$.
By continuity, there exists $\delta_1>0$ such that 
\bel{qj>}
q(t)>\hat q_1-\delta_1\,,\qquad q_j(t)>0\qquad\implies\qquad
\dot q_j(t)~\leq~-\ve_{1j}\,.\eeq
We observe that, if $q(t)>\hat q_1-\delta_1>0$, then $q_j(t)>0$ for some $j\in\O$.
Setting $\ve_1\doteq \min_j \ve_{1j}$, we obtain (\ref{impk}) for $\ell=1$.

{}From the implication
$$q(t)~\geq~\hat q_1-\delta_1\qquad\implies\qquad \dot q(t)~\leq~-\ve_1\,,$$
it follows 
$q(t)\leq \hat q_1-\delta_1$ for all $t\geq\tau_1$ sufficiently large. 
This yields (ii), for $\ell=1$.

Next, for $t>\tau_1$, the flux of cars arriving to the intersection 
from road 1 is
$$f_1(t)~=~\min\{ \omega_1^\diams\,,~c_1(M-q(t))\}$$
If road 1 is congested near the intersection, then for $t>\tau_1$ the outgoing flux is
$$f_1(t)~=~c_1(M-q(t))~\geq~
c_1(M-\hat q_1+\delta_1)\}~=~\omega_1^\diams-\delta_1'\,,$$
for some $\delta_1'>0$.
As a consequence,   road 1 must become free within time 
$$t_1~=~\tau_1+{1\over\delta_1'}\cdot\int_0^{\tau_1} [\omega_1^\diams - f_1(t)]\, dt\,.$$
This proves (iii), in the case $\ell=1$.
\v
The general inductive step is very similar.  Assume that the statements (i)--(iii)
have been proved for $\ell-1$.   Then for $t\geq t_{\ell-1}$ the incoming roads $i=1,\ldots,\ell-1$ are free. The flux of cars reaching the intersection
from these roads is $f_i(t)=\omega_i^\diams$.

Now assume that $t>t_{\ell-1}$ and $q(t)\geq \hat q_\ell$.    
In this case, 
 $q(t)\geq \hat q_i$ for all $i\in\I$, $i\geq \ell$.
Using Lemma 1,  for any $j\in \O$ we thus obtain
$$\dot q_j(t) ~\leq~\sum_{i< \ell} \omega_i^\diams \theta_{ij} + 
\sum_{i\geq \ell} c_i(M-q(t))\theta_{ij} - \omega_j^\diams\qquad\qquad \hbox{if}~~q_j(t)>0.$$
Therefore, if $q_j(t)>0$, then 
$$\dot q_j(t) ~\leq~-2\ve_{\ell j}~<~0$$
for some constants $\ve_{\ell j}$.
By continuity, there exists $\delta_\ell>0$ such that 
\bel{qjl>}
q(t)>\hat q_\ell-\delta_\ell\,,\qquad q_j(t)>0\qquad\implies\qquad
\dot q_j(t)~\leq~-\ve_{\ell j}\,.\eeq
Setting $\ve_\ell\doteq \min_j \ve_{\ell j}$, we obtain (\ref{impk}).

{}From the implication
$$q(t)~\geq~\hat q_\ell-\delta_\ell\qquad\implies\qquad \dot q(t)~\leq~-\ve_\ell\,,$$
it follows 
$q(t)\leq \hat q_\ell-\delta_\ell$ for all $t\geq\tau_\ell$ sufficiently large. 
This yields (ii).

Finally, for $t>\tau_\ell$, the flux of cars arriving to the intersection 
from road $\ell$ is
$$f_\ell(t)~=~\min\{ \omega_\ell^\diams\,,~c_\ell(M-q(t))\}$$
If road $\ell$ is congested near the intersection, then for $t>\tau_\ell$ 
the outgoing flux is
$$f_\ell(t)~=~c_\ell(M-q(t))~\geq~
c_\ell(M-\hat q_\ell+\delta_\ell)\}~=~\omega_\ell^\diams-\delta_\ell'\,,$$
for some $\delta_\ell'>0$.
As a consequence,  road $\ell$ must become free within time 
$$t_\ell~=~\tau_\ell+{1\over\delta_\ell'}\cdot\int_0^{\tau_\ell} [\omega_\ell^\diams - f_\ell(t)]\, dt\,.$$
This proves (iii).   By induction on $\ell$, our claim is proved.
\v
{\bf 2.} We now prove that, for any $\ve>0$,
there exists a time $t_\ve > t_\nu$ large enough so that
\bel{qub}
q(t)~\leq~q^*+\ve\qquad\qquad\forall t ~\geq~t_\ve\,.\eeq
Indeed, if $t>t_\ell$ , then the same arguments used before yield
the implication 
$$q(t)~\geq ~q^*+\ve\qquad\implies\qquad \dot q(t)~\leq~-\delta~<~0,$$
for some $\delta=\delta(\ve)>0$.
Hence, $q(t)\leq q^*+\ve$ whenever
$$t~\geq~t_\ve ~=~t_\nu +\delta^{-1}q(t_\nu) .$$ 
\v
{\bf 3.} In this step we prove a lower bound on the queue.
Namely, for any $\ve^\sharp>0$  there exists a time
$t^\sharp>0$ such that
\bel{qlb}
q(t)~\geq ~q^*-\ve^\sharp\qquad\qquad\forall t ~\geq~t^\sharp\,.\eeq
Toward this goal, choose  $j^*\in\O$ such that (\ref{gj*}) holds.
Then 
$$\dot q_{j^*}(t)~=~\sum_i\min\bigl\{ c_i(M-q(t)),\, 
\omega_i(\bar \rho(t))\bigr\}
\theta_{ij} -\bar f_{j^*}(t).$$
Now assume  $q(t)<q^*-\ve^\sharp$.   
If road $j^*$ is initially free, then it remains free for all times $t\geq 0$.
Hence $f_{j^*}(t)\leq f_j^{max}=\omega_j^\diams$.
In this case we have
$$\dot q_{j^*}(t)~\geq~\sum_i\min
\bigl\{ c_i(M-q^*+\ve^\sharp),~ \omega_i^\diams\bigr\}
\theta_{ij} -\omega_{j^*}^\diams~\geq~\delta^\sharp$$
for some $\delta^\sharp~>~0$. Therefore (\ref{qlb}) holds with
$$t^\sharp~\doteq~
{q_{j^*}(0)-q^*+\ve^\sharp\over\delta^\sharp}\,.$$

Next, we consider the more difficult case where the outgoing road $j^*$ is 
initially congested.
We claim that, if $q(\tau)\leq q^*$ at some time $\tau\geq t_\nu$,  then 
$q(t)\leq q^*$ for all $t\geq\tau$.
Indeed, for any $j\in\O$, if $q_j(t)>0$, then 
\bel{dq}\dot q_j(t)~\leq~\sum_{i\leq \nu} \omega_i^\diams \theta_{ij} + 
\sum_{i> \nu} c_i(M-q(t))\theta_{ij} - \omega_j^\diams\,.\eeq
Observing that the right hand side of 
(\ref{dq}) is nonpositive when $q(t)=q^*$, our claim is proved.

Now call
$$E_{j^*}(t)~\doteq~\omega_{j^*}^\diams \, t -\int_0^t \bar f_{j^*}(s)\, ds\,,$$
and observe that $E_{j^*}(t)\geq 0$ for all $t$, while
\bel{dde}
q(t)~\leq~q^*\qquad\implies\qquad \dot E_{j^*}(t)~\leq~0.\eeq
For $t>t_\nu$ we have 
\bel{qE}
\dot q_{j^*}(t) - \dot E_{j^*}(t)~\geq~\sum_{i\leq\nu} \omega_i^\diams\theta_{ij^*}
+\sum_{i>\nu}\min\left\{c_i(M-q(t)),~\omega_i^\diams\right\}\theta_{ij^*}
 - \omega_{j^*}^\diams\,.\eeq
If $q(t)\leq q^*-\ve^\sharp$, then
\bel{dqjj}\dot q_{j^*}(t) - \dot E_{j^*}(t)~\geq~\delta^\sharp\eeq
for some $\delta^\sharp>0$. 

Finally, assume that $q(\tau)<q^*-\ve^\sharp$, for some $\tau\geq t_\nu$.
Then  (\ref{dde}) and (\ref{dqjj}) imply (\ref{qlb}), with
$$t^\sharp~=~\tau + {E_{j^*}(\tau) + q^*-q(\tau)\over \delta^\sharp}\,.$$

\v
{\bf 4.} Denote by $\rho_k(t,x)$, $k\in \I\cup\O$,
the solution to the Riemann problem with buffer, and $\sigma_k(t,x)$
the self-similar solution determined by the Limit Riemann Solver (LRS).
{}From the limit $\lim_{t\to\infty}q(t)= q^*$ proved in the previous steps,
it follows that all boundary fluxes $ f_k(t)$ converge
to the corresponding boundary fluxes $\bar f_k$ in the 
self-similar solution determined by (LRS).  

Now consider  an incoming road $i\in \I$.
Since the initial data coincide
$$\rho_i(0,x)~=~\sigma_i(0,x)~=~\rho_i^\diams\qquad\qquad x<0\,,$$
for every $t>0$ by \cite{LF} we have the estimate
\bel{lfl}
\int_{-\infty}^0\bigl| \rho_i(t,x)- \sigma_i(t,x)\bigr|\, dx~\leq~
\int_0^t |f_i(s)- \bar f_i|\, ds\,.\eeq
 From the limit 
$$\lim_{t\to\infty} |f_i(t) - \bar f_i|~=~0$$
it thus follows
$$\lim_{t\to\infty} {1\over t} \int_{-\infty}^0\bigl| \rho_i(t,x)- \sigma_i(t,x)\bigr|\, dx~=~0.$$
For outgoing roads $j\in \O$, the estimates are  entirely similar.   
This achieves a proof of Theorem~2 in the case where (\ref{gj*}) holds for some $j^*\in\O$. 
\v
{\bf 5.} It remains to consider the case (Fig.~\ref{f:tf205}, left) where
\bel{efre}\sum_i \omega_i^\diams \theta_{ij}~<~\omega_j^\diams\eeq
for every $j\in \O$.  In this case, the arguments in step {\bf 1}
show that, for all $t\geq t_m$ sufficiently large, all incoming roads
become free.   In this case, for all times $t\geq t^\sharp$ sufficiently
large the incoming fluxes
are 
$$ f_i(t)~=~\omega_i^\diams~=~ \bar f_i\,.\qquad
\qquad i\in\I.$$   Moreover, for $t$ large all queue sizes
become $q_j(t)=0$, and the outgoing fluxes are 
$$ f_j(t)~=~\sum_{i}\omega_i^\diams\theta_{ij}~=~
\bar f_j\qquad\qquad j\in\O.$$
Inserting these identities in (\ref{lfl}), 
we conclude the proof as in the previous case.
\endproof


\v
{\bf Acknowledgment.} The first author  was partially supported
by NSF, with grant  DMS-1411786: ``Hyperbolic Conservation Laws and Applications". The second author recieved financial support for a stay at Penn State University from Tandberg radiofabrikks fond, Norges tekniske h\o gskoles fond, and Generaldirekt\o r Rolf \O stbyes stipendfond ved NTNU.
\vs

\end{document}